\documentclass[12pt,twoside]{amsart}
\usepackage{amsmath}
\usepackage{amssymb}
\usepackage{amscd}
\usepackage{xypic}
\xyoption{all}
\title[A Hodge index for Grothendieck residue pairing]{A Hodge index for Grothendieck residue pairing}

\author{Mohammad Reza Rahmati}
\thanks{}
\address{ Centro de Investigacion en Matematicas, A.C.
\hfill\break 
\hfill\break \\
\hfill\break }
\email{mrahmati@cimat.mx}

\newcommand{\comments}[1]{}

\newtheorem{theorem}{Theorem}[section]
\newtheorem{proposition}[theorem]{Proposition}
\newtheorem{corollary}[theorem]{Corollary}



\begin{document}

\begin{abstract}
In this text we apply the methods of Hodge theory for isolated hypersurface singularities to define a signature for the Grothendieck residue pairing of these singularities. 
\end{abstract}

\maketitle


\section*{Introduction}

A Hodge structure (HS) of weight $n \in \mathbb{Z}$ on a $\mathbb{Q}$-vector space $H$ consists of a decomposition $H \otimes \mathbb{C}=\oplus_{p+q=n}H^{p,q}$ into $\mathbb{C}$-subspaces $H^{p,q}$ such that $H^{q,p}=\overline{H^{p,q}}$. This structure can equally been explained by the existence of a decreasing flag $F^{\bullet} : F^0=H \supset F^1 \supset ... \supset F^n=0$ such that $F^i \cap \overline{F^{n-i+1}}=0$. The HS is said to be polarized if there exists a non-degenerate $\mathbb{Q}$-bilinear form $S:H \otimes H \to \mathbb{Q}$ such that the aforementioned decomposition is orthogonal with respect to this bilinear form and the the bilinear form $S(.,C.)$ is positive definite, where $C:H \to H$ defined by $C\mid_{H^{p,q}}=i^{p-q}$ is the Weil operator. A mixed Hodge structure (MHS) is given by a triple $(H,F^{\bullet},W_{\bullet})$ where $H$ and $F^{\bullet}$ are as before and $W_{\bullet}$ is an increasing filtration on $H$ defined over $\mathbb{Q}$. This data is subject to the condition that for all $k$, the pairs $(Gr_k^W H, F^{\bullet}Gr_k^W)$ is a HS of weight $k$. The polarization for a MHS is given by a set polarizations $S_k$ on the pure HS $(Gr_k^W H, F^{\bullet}Gr_k^W)$ for each $k$. Mixed HS's constitute a category, with morphisms to be linear maps respecting the filtrations. Examples of HS (resp. MHS) are cohomologies of projective (resp. quasi-projective) varieties, \cite{SCH}, \cite{D1}. A well known theorem of P. Deligne states that any MHS $H$ has a unique functorial bigrading $H=\bigoplus_{p,q}I^{p,q}$ such that

\begin{equation}
\overline{I^{p,q}}=I^{q,p} \qquad \mod \bigoplus_{r<q,s<p}I^{r,s}
\end{equation}

We propose to study families of HS of the same weight over the puncture disc. An example of this happens when one considers a projective map $f: X \to S$ between quasi-projective varieties. 
In this case the $k$-th cohomology of each fiber has a polarized HS of weight $k$. The special case we will consider in this text is when $f: \mathbb{C}^{n+1} \to \mathbb{C}$ has isolated singularity and is given by a holomorphic germ around $0 \in \mathbb{C}^n$. By a theorem of J. Milnor $f$ has a Milnor fibration representative $f:X \to \Delta^*$ defined on sufficiently small neighborhoods of zero. The only interesting non-trivial cohomology of the fibers, $X_t:=f^{-1}(t)$ (in this case) is the middle cohomology, i.e. $n$-th cohomology (also called vanishing cohomology). Let $\exp:U \to \Delta^*$ is the universal covering of the punctured disc, it is convenient to set $X_{\infty}:=X_t \times_{\Delta^*} U$ called canonical fiber. 

The subject of this article is to study the variation of MHS obtained from the vanishing cohomology of such fibrations, and try to relate them to the Grothendieck residue pairing associated to $f$. 
The hypothesis that $f$ has an isolated singularity implies that the module $\Omega_f=\Omega_X^{n+1}/df \wedge \Omega_X^n \cong \mathcal{O}_X/\partial f$ has finite rank equal to the Milnor number of $f$. The Grothendieck pairing (local residue) is a symmetric $\mathbb{C}$-bilinear pairing 

\begin{equation}
res_{f,0}: \Omega _{f}\times \Omega _{f}\to \mathbb{C}
\end{equation}

It is a well known theorem of A. Grothendieck that the local residue form is non-degenerate. A natural question is; how to calculate the signature of this form by the topological invariants of the singularity. This is the subject of this note and we try to relate this  to the Hodge numbers of the MHS defined by $f$. We explain that the $\mathbb{C}$-vector space $\Omega_f$ inherits a MHS from $H^n(X_{\infty},\mathbb{C})$ by constructing a specific isomorphism between them. Then we show that a little modification of the local residue (2) can be regarded as a polarization for the MHS.

\section{Variation of mixed Hodge structure}

A polarized variation of mixed Hodge structure over the punctured disc $\Delta^*$ consists of the 5-tuple $( \mathcal{H}, F^{\bullet}, W_{\bullet}, \nabla, S)$ where 

\begin{itemize}

\item $\mathcal{H}$ is local system of $\mathbb{Q}$-vector spaces on $\Delta^*$.
\item $W_{\bullet}=(W_i)$ is an increasing filtration on $\mathcal{H}$ by sub-local systems of $\mathbb{Q}$-vector spaces.
\item $F^{\bullet}=(F^i)$ is a decreasing filtration on the vector bundle $\mathcal{H} \otimes_{\mathbb{Q}} \mathcal{O}_{\Delta^*}$ by holomorphic sub-bundles.
\item $\nabla: \mathcal{H} \otimes_{\mathbb{Q}} \mathcal{O}_{\Delta^*} \to \mathcal{H} \otimes_{\mathbb{Q}} \Omega_{\Delta^*}^1$ is a flat connection satisfying Griffiths transversality;

\begin{equation}
\nabla(F^i) \subset F^{i-1} \otimes \Omega_{\Delta^*}^1
\end{equation}

\item $S: \mathcal{H} \times \mathcal{H} \to \mathbb{Q}$ is a flat pairing inducing a set of rational non-degenerate bilinear forms $S_k:Gr_W^k \mathcal{H} \otimes Gr_W^k \mathcal{H} \to \mathbb{Q}$ such that the triple $(Gr_W^k \mathcal{H}, F^{\bullet}Gr_W^k ,S_k)$ defines pure polarized variation of Hodge structure on $\Delta^*$. 

\end{itemize}

\noindent
which we briefly mention as $\mathcal{H}$. In case explained the flat connection $\nabla$ is usually denoted by $\partial_t$.

If $f:(\mathbb{C}^{n+1},0) \to (\mathbb{C},0)$ is a holomorphic germ having an isolated singularity at $0 \in \mathbb{C}^{n+1}$, then the family of middle cohomologies of the fibers $X_t:=f^{-1}(t)$ constitute a local system $\mathcal{H}:= R^nf_* \mathbb{Q}=\bigcup_t H^n(X_t,\mathbb{Q})$, of rank $\mu(f)$ the Milnor number of $f$. The monodromy $M:H^n(X_t,\mathbb{Q}) \to H^n(X_t,\mathbb{Q})$ around $0 \in \mathbb{C}$ can be written as $M=M_sM_u$ where $M_s$ is a diagonal matrix called the semisimple part of $M$. The logarithm $N=\log M_u \otimes 1/2 \pi i$ is a nilpotent map on $H_t=H^n(X_t,\mathbb{Q})$. By a well known theorem of Jacobson-Morosov \cite{H1}, $N$ induces an increasing filtration $W_{\bullet,t}$ on $H_t$. We will use the notation $X_{\infty}=X_t \times_{\Delta^*} U$ where $U$ is the upper half plane. This changes nothing in the middle cohomologies. 

According to the Riemann-Hilbert correspondence, the local system $\mathcal{H}$ corresponds to Gauss-Manin system $(\mathcal{G}=\mathcal{H} \otimes \mathcal{O}_{\Delta^*}, \nabla)$ on the punctured disc, with a flat connection

\begin{equation}
\partial_t:= \nabla_{d/dt}:\mathcal{H} \to \mathcal{H}
\end{equation}
 
\begin{proposition} \cite{SCHU}, \cite{SC2}
If $f:(\mathbb{C}^{n+1},0) \to (\mathbb{C},0)$ is a holomorphic germ having an isolated singularity at $0 \in \mathbb{C}^{n+1}$, then the action of $\partial_t$ in (4) is invertible.
\end{proposition}

We will briefly explain the methodology to understand this theorem in what follows. The $\mathbb{C}[t,t^{-1}]$-module $\mathcal{G}$ has a canonical filtration due to Malgrange-Kashiwara, namely $V$-filtration indexed by $\alpha \in \mathbb{Q}$. It is characterized by the properties; 

\begin{itemize}
\item $t.V^{\alpha} \subset V^{\alpha+1}$, 
\item $\partial_t .V^{\alpha} \subset V^{\alpha-1}$ 
\item the operator $t\partial_t-\alpha$ is nilpotent on $Gr_V^{\alpha}$. 
\end{itemize}

\noindent
The $V$-filtration always exists and is unique and its definition is independent of choice of the coordinate $t$, and we set $C^\alpha=V^\alpha/V^{>\alpha}$. The Gauss-Manin system can also be written as;

\begin{equation}
\mathcal{G}=\displaystyle{\sum_{-1 \leq \alpha <0}}\mathbb{C}\{t\}[t^{-1}]C^{\alpha}.
\end{equation}

\noindent
There is an isomorphism 

\begin{equation}
\psi := \displaystyle{\bigoplus_{-1<\alpha \leq 0}}\psi^\alpha :\displaystyle{\bigoplus_{-1<\alpha \leq 0}}H_{\mathbb{C}}^{\lambda_\alpha} \to \displaystyle{\bigoplus_{-1<\alpha \leq 0}}C^\alpha, \qquad \psi_{\alpha}(A_\alpha)= t^{\alpha}\exp(\log t.{-N \over 2 \pi i}).A_{\alpha}
\end{equation}

\noindent
where $H_{\mathbb{C}}^{\lambda}:=H^n(X_{\infty},\mathbb{C})^{\lambda}=\ker(M_s-\lambda)^{r} \subset H^n(X_{\infty},\mathbb{C})$ is the generalized eigenspaces of the monodromy around $0 \in \Delta^*$. We have $\lambda_{\alpha}:=\exp (2 \pi i \alpha), \ -1 < \alpha \leq 0$. The action of the monodromy operator $M$ on $H_{\mathbb{C}}^{\lambda_{\alpha}}$ correspond to the action of the operator $ \exp(t\partial_t-\alpha)$ by this map. The explanation of the $V$-filtration is simple in this case and $V^{\alpha}$ is the submodule generated by the image of sections with index not less than $\alpha$. 

E. Brieskorn \cite{B} considers the $\mathcal{O}_{\Delta^*}$-modules
 
\begin{equation}
H^{\prime \prime}= \displaystyle{\frac{f_*\Omega_X^{n+1}}{df \wedge d(\Omega_X^{n-1})}}
\end{equation}

\noindent
called the Brieskorn lattice associated to the fibration of $f$. Sometimes the stack of the sheaf $H''$ is called the Brieskorn lattice  of $f$. It is a locally free sheaf of rank $\mu(f)$ the Milnor number of $f$. There is a map relating the sections of the sheaf $\Omega_X$ to sections of $\mathcal{G}$ namely 

\begin{equation}
s:\omega \longmapsto \frac{\omega}{df}=\{ t \mapsto Res_{(f=t)} \frac{\omega}{f-t} \} 
\end{equation}

\noindent
which takes values in $V^{>-1}$, and its kernel is $df \wedge d\Omega_X^{n-1}$. Therefore it embeds $H''$ into $V^{>-1} \subset \mathcal{G}$. Using the $V$-filtration on the Gauss-Manin module; we may define the Hodge filtration on $H^n(X_{\infty},\mathbb{C})$ by

\begin{equation}
F^pH^n(X_{\infty},\mathbb{C})_{\lambda}=\psi_{\alpha}^{-1}\displaystyle{(\frac{V^{\alpha} \cap \partial_t^{n-p}H_0^{\prime \prime}}{V^{>\alpha}})}, \qquad \alpha \in (-1,0] 
\end{equation}

where the action of $\partial_t$ is via the embedding induced by $s$ and $F^pH^n(X_{\infty}) =\oplus_{\lambda} F^pH^n(X_{\infty},\mathbb{C})_{\lambda}$, \cite{H1}, \cite{R}, \cite{SCHU}, \cite{KUL}, \cite{SC2}. 

Finally, in order to define the polarization form we proceed as follows. There exists a change of variable such that the holomorphic map $f$ becomes a polynomial of sufficiently high degree $d$. Then, the Milnor fibration $f:X \to \Delta^*$ can be embeded into a projective fibration $f_Y:Y \to \Delta^*$ such that Zero is the only singular point of the closure $Y_{0}=\overline{X_0}^{\text{Zar}}$ in a hyperplane in $P^{n+1}(\mathbb{C})$ and the fibers in $Y$ are prolongations of the fibers in $X$. The map induced on the middle cohomologies from the embedding $i:X_{\infty} \hookrightarrow Y_{\infty}$ is surjective and we have the exact sequence

\begin{equation}
0 \to \ker(M_Y-id) \to P^n(Y_\infty) \to H^n(X_\infty) \to 0.
\end{equation}

\noindent
where $P^n(Y_{\infty})$ is the primitive part of $H^n(Y_{\infty})$ and (10) is a form of the invariant cycle theorem. The monodromy $M_Y$ and its logarithm $N_Y$ are similarly defined for the fibration $Y$. The polarization form is given by $S=S_{ \ne 1} \oplus S_1$ where

\begin{center}
$ S_{\ne 1}(a,b) = S_Y(i^*a,i^*b), \qquad a,b \in H^n(X_{\infty})_{\ne 1} $ \\[0.2cm]
$S_{1}(a,b) = S_Y(i^*a,i^*N_Yb), \qquad a,b \in H^n(X_{\infty})_{1} $
\end{center}

\noindent
where $S_Y$ is the polarization form for the HS $H^n(Y_{\infty})$. The first two vector spaces in (10) have a natural MHS. J. Steenbrink establishes that there exists a unique MHS on $H^n(X_{\infty})$ which makes the short sequence to be exact in the category of MHS's, called it Steenbrink limit MHS, \cite{JS2}, \cite{SCH}, \cite{SC2}, \cite{H1}, \cite{JS3}.

\begin{theorem} \cite{H1}
The filtration defined in (9) is the Steenbrink limit Hodge filtration. The 5-tuple $(\mathcal{H},F^{\bullet},W_{\bullet},\partial_t,S)$ constitute a variation of mixed Hodge structure on the punctured disc.
\end{theorem}

\section{Grothendieck residue pairing}

For a holomorphic germ $f:(\mathbb{C}^{n+1},0) \to (\mathbb{C},0)$ with an isolated critical point; the local (Grothendieck) residue is 
 
\begin{center}
$g \longmapsto Res_{0}{\displaystyle{\left[\frac{gdx}{\frac{\partial f}{\partial{x_0}}...\frac{\partial{f}}{\partial{x_n}}}\right]}}:= \displaystyle{\frac{1}{(2\pi {i})^{n+1}}} \displaystyle{\int _{\Gamma_{\varepsilon}}} \frac{gdx}{\frac{\partial{f}}{\partial{x_0}}...\frac{\partial{f}}{\partial{x_n}}}$
\end{center}

\noindent
where $\Gamma_{\epsilon}$ is a suitable cocycle around $0 \in \mathbb{C}^{n+1}$. One can define a bilinear form $res_{f,0}$
 
\begin{equation}
res_{f,0}: \Omega _{f}\times \Omega _{f}\to \mathbb{C}
\end{equation}
\begin{center}
$({g_{1}}dx,{g_{2}}dx) \longmapsto Res_{0}{\displaystyle{\left[\frac{g_{1}g_{2}dx}{\frac{\partial f}{\partial{x_0}}...\frac{\partial{f}}{\partial{x_n}}}\right]}}$,
\end{center}

\noindent
It is a symmetric bilinear pairing (Grothendieck residue pairing-local residue), which is non-degenerate, \cite{G3}. There is a simple relation between the module of relative differentials $\Omega_f$ and the Brieskorn lattice $H''$ given by

\begin{equation}
\frac{H''}{\partial_t^{-1}H''} \cong \Omega_f
\end{equation}

In \cite{R} we have defined an isomorphism $\Phi:H^n(X_{\infty}) \to \Omega_f$ which allows to transfer the Hodge and weight filtration over $\Omega_f$. According to the theorem of Deligne-Hodge the MHS $H^n(X_{\infty})$ has a decomposition $H^n(X_{\infty})=\bigoplus I^{p,q}$ satisfying (1). The map $\Phi$ is defined as follows on the generalized eigenspace in $I_{\lambda}^{p,q}$.

\begin{equation}
\Phi_{\lambda}^{p,q}: I^{p,q}_{\lambda} \stackrel{\hat{\Phi}_{\lambda}}{\longrightarrow}Gr_V^{\alpha+n-p}H'' \stackrel{pr}{\longrightarrow} Gr_V^{\bullet} H^{\prime \prime}/\partial_t^{-1}H^{\prime \prime} \stackrel{\cong}{\longrightarrow} \Omega_f
\end{equation}

\noindent
where $\hat{\Phi}_{\lambda}^{p,q}:= \partial_t^{p-n} \circ \psi_{\alpha}$. The last isomorphism forgets the $V$-grading on $\Omega_f$, and in this way $\Phi$ is not unique if defined. A choice of $\Phi$ corresponds to a choice of a section of 

\begin{equation} 
Gr_V^{\alpha}[V^{\alpha} \cap H''] \to Gr_V^{\alpha}[H''/\partial_t^{-1}H''] 
\end{equation}

\noindent
for $-1 \leq \alpha <0$, \cite{SAI6}, \cite{R}, \cite{H1}. If $J^{p,q}=\Phi^{-1} I^{p,q}$ is the corresponding subspace of $\Omega_f$, then $\Omega_f=\bigoplus_{p,q}J^{p,q}$ and define $\tilde{C}|_{J^{p,q}}:=(-1)^{p}$. The following diagram is commutative (\cite{R} Theorem 8.6.1);

\begin{equation}
\begin{CD}
\widehat{Res}_{f,0}:\Omega_f \times \Omega_f @>>> \mathbb{C}\\
\ \ \ \ \ @VV(\Phi^{-1},\Phi^{-1})V                   @VV \times *V \\
S:H^n(X_{\infty}) \times H^n(X_{\infty}) @>>> \mathbb{C}
\end{CD} \qquad \qquad  * \ne 0
\end{equation}

\noindent
where, 

\[ \widehat{Res}_{f,0}=\text{res}_{f,0}\ (\bullet,\tilde{C}\ \bullet) \]

\noindent
In other words;

\begin{equation}
S(\Phi^{-1}(\omega),\Phi^{-1}(\eta))= * \times \ \text{res}_{f,0}(\omega,\tilde{C}.\eta), \qquad 0 \ne * \in \mathbb{C}.
\end{equation}

The mixed Hodge structure on $\Omega_f$ is defined via $\Phi$, \cite{R} (see also \cite{V} and \cite{SAI6}). In case $f=f(z)$ is a quasihomogeneous polynomial, the mixed Hodge structure can be explained via monomial basis of the Milnor or Jacobi ring. In this case the inverse map $\Phi^{-1}$ is as follows,

\begin{equation}
\Phi^{-1}: [z^{\alpha} dz] \longmapsto c_{\alpha}. [res_{f=1}( z^{\alpha} dz/(f-1)^{[l(\alpha)]}) ]
\end{equation}

\noindent
with $c_{\alpha} \in \mathbb{C}$, and $z^{\alpha}$ is a monomial basis of Jacobi ring of $f$ by the work of Griffiths and Steenbrink. In this case the Hodge structure is the same as pole filtration can be explained by the degrees $l(\alpha)=\sum (\alpha_i+1)w_i$ where $w_i$ is the weight of $z_i$, \cite{JS7}, \cite{R}.

\begin{proposition} \cite{R}(Riemann-Hodge bilinear relations for Grothendieck residue pairing on $\Omega_f$) Assume $f:\mathbb{C}^{n+1} \to \mathbb{C}$ is a holomorphic germ with an  isolated singularity. The 3-tuple $(\Omega_f, \Phi F^{\bullet},\Phi W_{\bullet})$ define a polarized MHS, which is polarized by $\widehat{Res}_{f,0})$ in the following sense. Suppose $\mathfrak{f}$ is the corresponding map to $N:H^n(X_{\infty}) \to H^n(X_{\infty})$, via the isomorphism $\Phi$. Define 

\[ P_l=PGr_l^W:=\ker(\mathfrak{f}^{l+1}:Gr_l^W\Omega_f \to Gr_{-l-2}^W\Omega_f) \]

\noindent
Going to $W$-graded pieces;
\begin{equation}
\widehat{Res}_l: PGr_l^W \Omega_f \otimes_{\mathbb{C}} PGr_l^W \Omega_f \to \mathbb{C}
\end{equation}

\noindent
is non-degenerate and according to Lefschetz decomposition 

\[ Gr_l^W\Omega_f=\bigoplus_r \mathfrak{f}^r P_{l-2r} \]

\noindent
we will obtain a set of non-degenerate bilinear forms,

\begin{equation}
\widehat{Res}_l \circ (id \otimes \mathfrak{f}^l): P Gr_l^W \Omega_f  \otimes_{\mathbb{C}} P Gr_l^W \Omega_f  \to \mathbb{C}, 
\end{equation} 
\begin{equation}
\widehat{Res}_l=res_{f,0}\ (id \otimes \tilde{C} .\  \mathfrak{f}^l)
\end{equation}

\noindent
such that the corresponding hermitian form associated to these bilinear forms is positive definite. In other words, 

\begin{itemize}

\item $\widehat{Res}_l(x,y)=0, \qquad x \in P_r, \ y  \in P_s, \ r \ne s $

\item If $x \ne 0$ in $P_l$, 

\[ Const \times res_{f,0}\ (C_lx,\tilde{C} .\  \mathfrak{f}^l .\bar{x})>0  \]

\noindent
where $C_l$ is the corresponding Weil operator. 

\end{itemize}

\end{proposition}

Hodge theory assigns to any polarized Hodge structure $(H,F,S)$ a signature which is the signature of the hermitian form $S(C \bullet, \bar{\bullet})$, where $C$ is the Weil operator. In case of a polarized mixed Hodge structure $(H,F,W(N),S)$, where $N$ is a nilpotent operator this signature is defined to be the sum of the signatures of the hermitian forms associated to the graded polarizations $S_l: PGr_l^W H \times PGr_l^W H \to \mathbb{C}$, i.e signatures of $h_l:=S_l(C_l \bullet , N^l \bar{\bullet})$ for all $l$. The signature associated to the polarized variation of mixed Hodge structure of an isolated hypersurface singularity with even dimensional fibers is calculated in \cite{JS2} as 

\begin{equation} 
\sigma=\displaystyle{\sum_{p+q=n+2}(-1)^q h_1^{pq}+2\sum_{p+q \geq n+3}(-1)^q h_1^{pq}+\sum_{p,q} (-1)^q h_{\ne 1}^{pq}} 
\end{equation}

\noindent
where $h_1=\dim H^n(X_{\infty})_1,h_{\ne 1}=\dim H^n(X_{\infty})_{\ne 1}$ are the corresponding Hodge numbers. This signature is $0$ when the fibers have odd dimensions.

\begin{corollary}  The signature associated to the Grothendieck pairing $\ {res}_{f,0}\ $ of an isolated hypersurface singularity germ $f$; is equal to the signature of the polarization form associated to the MHS of the vanishing cohomology and is given by (21). This index is zero when fibers have odd dimensions.
\end{corollary}

\begin{proof}
This follows from the diagram (15) and Theorem 1.2 and proposition 2.1.
\end{proof}

In the Hodge terminology this index can also be associated to the real hypersurface $f:\mathbb{R}^{n+1} \to \mathbb{R}$ when the isolated singularity is algebraic, i.e. it remains isolated over $\mathbb{C}$.

\end{document}